\documentclass[10pt]{amsart}

\usepackage{amsmath,amsfonts,amssymb,amsthm,amscd,epsfig}
\usepackage{youngtab,comment}

\def\C{\mathbb{C}}
\def\CP{\mathbb{CP}}
\def\H{\mathbb{H}}
\def\Z{\mathbb{Z}}

\def\T{T^2}
\def\id{\mathbf{1}}

\def\v{\mathbf{v}}
\def\tv{{\tilde{v}}}
\def\w{\mathbf{w}}
\def\e{\mathbf{e}}

\def\B{\mathbf{B}}

\def\M{\mathfrak{M}}

\def\B{\mathfrak{B}}

\def\p{\mathfrak{p}}
\def\P{\mathcal{P}}

\def\Cl{\mathrm{Cl}}

\newcommand{\vac}[1]{ \left| #1 \right> }
\newcommand{\hs}[1][n]{X_{#1}}

\DeclareMathOperator{\End}{End}

\DeclareMathOperator{\tr}{tr}

\DeclareMathOperator{\Id}{Id}

\newtheorem{theo}{Theorem}[section]
\newtheorem{prop}[theo]{Proposition}

\newtheorem*{rem*}{Remark}

\numberwithin{equation}{section}

\begin{document}
\title{A geometric boson-fermion correspondence}
\author{Alistair Savage}
\address{University of Ottawa\\
Ottawa, Ontario \\ Canada} \email{alistair.savage@uottawa.ca}
\thanks{This research was supported in part by the Natural
Sciences and Engineering Research Council (NSERC) of Canada}
\subjclass[2000]{Primary: 14C05,17B69; Secondary: 55N91}
\keywords{Boson-fermion correspondence, vertex algebras, affine Lie
algebras, quivers, quiver varieties, Hilbert schemes, equivariant
cohomology}

\begin{abstract}
The fixed points of a natural torus action on the Hilbert schemes of
points in $\C^2$ are quiver varieties of type $A_\infty$.  The
equivariant cohomology of the Hilbert schemes and quiver varieties
can be given the structure of bosonic and fermionic Fock spaces
respectively.  Then the localization theorem, which relates the
equivariant cohomology of a space with that of its fixed point set,
yields a geometric realization of the important boson-fermion
correspondence.

\vspace{0.5cm}

\paragraph{\bf Version fran\c{c}aise:} Les points fixes d'une action
canonique d'un tore sur le sch\'ema de Hilbert de $\C^2$ sont des
vari\'et\'es de quiver de type $A_\infty$.  On peut donner la
cohomologie \'equivariante des sch\'emas de Hilbert et des
vari\'et\'es de quiver la structure des \'espaces de Fock
fermionique et bosonique respectivement.  Alors, la th\'eor\`eme de
localisation, qui lie la cohomologie \'equivariante d'une \'espace
avec la cohomologie \'equivariante de son ensemble des point fixes,
nous permet de donner une r\'ealisation g\'eom\'etrique de la
correspondance bosonique-fermionique.
\end{abstract}

\maketitle

\section*{Introduction}
Recently there has been substantial interest in geometric
constructions in representation theory. Such constructions translate
between purely algebraic representation theoretic statements and
statements involving geometric objects such as flag varieties,
affine Grassmannians, quiver varieties and Hilbert schemes. This
often provides one with new geometric techniques to examine various
topics in representation theory (such as in the proof of the
Kazhdan-Lusztig conjecture) as well as representation theoretic
tools to organize and study the structure of various geometric
objects.

In this paper, we will focus on two particular geometric
constructions.  One of these involves varieties associated to
quivers.  Quivers, which are simply directed graphs, and their
representations have a long history (see \cite{Sav06a}).  Lusztig
\cite{L91} associated certain varieties to quivers and used these to
provide a geometric realization of half of the universal enveloping
algebra (or its quantum analogue) of Kac-Moody algebras. Then in
\cite{N94,N98}, Nakajima modified these quiver varieties and gave a
geometric construction of the representations of these algebras. The
underlying vector space of the representation is the homology of the
quiver varieties and the action of the Kac-Moody algebra is given by
certain \emph{correspondences} in products of these varieties.
Nakajima's construction was motivated by his work with Kronheimer on
solutions to the anti-self-dual Yang-Mills equations on ALE
gravitational instantons \cite{KN90}.

The second geometric construction we consider in this paper was
developed by Nakajima \cite{Nak99} and Grojnowski \cite{Gro96}.  It
realizes irreducible representations of infinite-dimensional
Heisenberg algebras in the (co)homology of the Hilbert schemes of
points on surfaces.  This is very similar to the quiver variety
picture. In fact, one can view the Hilbert schemes as quiver
varieties associated to the so-called Jordan quiver.

\[
\includegraphics[height=0.15\textwidth]{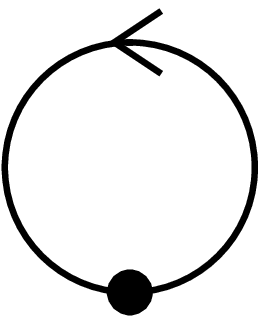}
\]

These geometric realizations have proven to be very useful. In
particular, they provide us with remarkable bases for the
corresponding algebraic objects (the so-called canonical and
semicanonical bases) that have very nice positivity, integrality and
compatibility properties.  They also have beautiful connections to
other areas of mathematics such as the theory of crystals.  Their
discovery has opened up a fruitful line of research -- to extend the
geometric viewpoint by developing geometric constructions of other
algebraic results.  This approach has been successful in several
areas including Weyl group actions, representations of the Virasoro
algebra, Demazure modules, and Clifford algebras.  In the current
paper, we continue along this path and provide a geometric
construction of the so called \emph{boson-fermion correspondence}.

The boson-fermion correspondence is of fundamental importance in
mathematical physics. In physics, the terms \emph{boson} and
\emph{fermion} refer to particles of integer and half-integer spin
respectively.  \emph{Bosonic} and \emph{fermionic Fock space}, which
can be thought of as certain state spaces of bosons and fermions,
are mathematically defined to be representations of an
infinite-dimensional Heisenberg or oscillator algebra and an
infinite-dimensional Clifford algebra (and the Lie algebra
$sl_\infty$ or $gl_\infty$) respectively.  The boson-fermion
correspondence describes a precise relationship between these two
spaces.  It plays an important role in the theory of vertex
operators and the basic representation of affine Lie algebras (see
\cite{DJKM82,Fre81} and references therein).

We will use the constructions discussed above to describe a
realization of the boson-fermion correspondence using the geometry
of Hilbert schemes and quiver varieties.  Let $\hs$ denote the
Hilbert scheme of $n$ points in $\C^2$ (see
Section~\ref{sec:hilbert-scheme-fixed-pts}). There is a natural
torus action on $\hs$ and one can consider the associated
equivariant cohomology on which there is an action of an
infinite-dimensional Heisenberg algebra. The generators of the
algebra act by ``adding or removing points'' along the $x$-axis in
$\C^2$. This yields the geometric realization of bosonic Fock space.
That of the fermionic Fock space is obtained by considering the
quiver varieties corresponding to the basic representation of
$sl_\infty$. In this case the varieties are simply points which, as
was shown in \cite{FS03}, can be naturally enumerated by Young
diagrams.

We will see that the torus fixed points of the Hilbert schemes are
naturally identified with the $sl_\infty$ quiver varieties. The
localization theorem states that under certain assumptions the
equivariant cohomology of a space with the action of a torus $T$ is
isomorphic to the equivariant cohomology of its $T$-fixed points.
Thus, the localization theorem yields an isomorphism between the
geometric constructions of the bosonic and fermionic Fock spaces and
a geometric boson-fermion correspondence. We then see that the
bosons correspond to ``global'' additions of points while the
fermions corresponds to ``local'' operators.  That is, the bosons
naturally act on the equivariant cohomology of the entire Hilbert
scheme while the fermions naturally act at its torus fixed points.

We note that in \cite{LQW04} Li, Qin, and Wang found that the
\emph{multi-point trace function}, a generating function of
intersection numbers of equivariant Chern characters in spaces
isomorphic to the Hilbert schemes mentioned above, is related in a
simple way to the characters of the fermionic Fock space when the
equivariant cohomology of these spaces (with the natural geometric
structure of bosonic Fock space) is identified with fermionic Fock
space via the boson-fermion correspondence.

As mentioned above, the geometric constructions of the two
representations considered in this paper are not new.  The goal here
is rather to examine the interplay between them.  The important step
is the identification of the quiver varieties for $sl_\infty$ with
the torus-fixed points of the Hilbert scheme.  We believe this point
of view to be important for two reasons.  First of all, it allows
one to identify a fundamental concept in equivariant cohomology, the
localization theorem, with an important concept in mathematical
physics, the boson-fermion correspondence.  Secondly, we expect this
idea to lead to other interesting results.  For example, if one
considers a finite cyclic subgroup $\Gamma$ of the torus, the
$\Gamma$-fixed points of the Hilbert scheme can be naturally
identified with affine quiver varieties of type $A$. Thus, we expect
that an extension of the ideas presented here to this setting will
lead to new geometric interpretations of the vertex operator
construction of representations of affine Lie algebras (see
\cite{Gro96} for some results in this direction) and perhaps
explicit algebraic descriptions in the vertex operator framework of
the natural bases coming from the geometric picture. We see the
current paper as an important first step in this direction.

The organization of the paper is as follows.  In
Sections~\ref{sec:boson-fermion} and~\ref{sec:localization} we
review the boson-fermion correspondence and the localization theorem
in equivariant cohomology.  In
Section~\ref{sec:geometric-fock-space} we describe the geometric
construction of fermionic Fock space using quiver varieties and in
Section~\ref{sec:hilbert-scheme-fixed-pts} we describe the torus
fixed points of the Hilbert scheme and their identification with
quiver varieties.  Then in Section~\ref{sec:geom-bosonic} we recall
the geometric realization of bosonic Fock space using Hilbert
schemes. Finally, in Section~\ref{sec:geom-correspondence} we state
the geometric boson-fermion correspondence.

The author would like to thank I. B. Frenkel for many invaluable
comments and suggestions.  He would also like to thank M. Harada, A.
Knutson, M. Haiman, and W. Wang for answering his questions on
equivariant cohomology and M. Lau for useful discussions.  Finally,
he would like to thank the American Institute of Mathematics and the
organizers of the ARCC Workshop on Generalized Kostka Polynomials
where some of the ideas in this paper were developed.


\section{The Boson-Fermion Correspondence}
\label{sec:boson-fermion}

The \emph{boson-fermion correspondence} is an isomorphism between
two representations of an infinite-dimensional Heisenberg algebra
(or the closely related \emph{oscillator algebra}). These
representations are on the bosonic and fermionic Fock spaces. In
this section we recall the two representations and the isomorphism
between them.  For further details we refer the reader to
\cite[{\S\S\ 14.9-14.10}]{K}.

Define the \emph{oscillator algebra} to be the Lie algebra
\[
\mathfrak{s} = \bigoplus_{m \in \Z} \C s_m \oplus \C K
\]
with commutation relations
\begin{equation} \label{s-relations}
[\mathfrak{s},K]=0, \quad [s_m,s_n] = m \delta_{m,-n} K.
\end{equation}
The subalgebra spanned by $s_n$, $n \ne 0$, and $K$ is an
infinite-dimensional Heisenberg algebra.  The oscillator algebra has
a natural representation on the \emph{full bosonic Fock space}
\[
B = \C[p_1,p_2,\dots;q,q^{-1}],
\]
a polynomial algebra on indeterminates $p_1, p_2, \dots$ and
$q,q^{-1}$.  Physically, this space can be thought of as a certain
state space of bosons (particles of integer spin).  The
indeterminate $p_k$ represents a particle in state $k$.  Note that
more than one boson can occupy the same state.

The representation $r^B$ of the oscillator algebra $\mathfrak{s}$ on
$B$ is given by:
\begin{gather*}
r^B(s_m) = m\frac{\partial}{\partial p_m},\quad r^B(s_{-m}) = p_m
\text{ for } m>0, \\
r^B(s_0) = q \frac{\partial}{\partial q},\quad r^B(K) = 1.
\end{gather*}
The operators $s_{-m}$ and $s_m$ can be thought of as creation and
annihilation operators respectively.

We next describe another representation of the oscillator algebra.
An infinite expression of the form
\[
\underline{i_0} \wedge \underline{i_1} \wedge \underline{i_2} \wedge
\dots,
\]
where $i_0, i_1, \dots$ are integers satisfying
\[
i_0 > i_1 > i_2 > \dots,\ i_n = i_{n-1} - 1 \text{ for } n \gg 0,
\]
is called a \emph{semi-infinite monomial}.  Let $F$ be the complex
vector space with basis consisting of all semi-infinite monomials,
and let $H(\cdot,\cdot)$ denote the Hermitian form on $F$ for which
this basis is orthonormal.  $F$ is called \emph{full fermionic Fock
space}.  Physically, this can be thought of as a certain state space
of fermions (particles of half-integer spin) with the integers
appearing in a semi-infinite monomial labeling the various states.
Note that the fermions satisfy the \emph{Pauli exclusion principle}
-- no two particles can occupy the same state. Let
\[
\vac{m} = \underline{m} \wedge \underline{m-1} \wedge \underline{m-2} \wedge \dots
\]
be the \emph{vacuum vector} of charge $m$.  We say that a
semi-infinite monomial has \emph{charge} $m$ if it differs from
$\vac{m}$ at only a finite number of places.  Thus $\varphi =
\underline{i_0} \wedge \underline{i_1} \wedge \dots$ is of charge
$m$ if $i_k = m - k$ for $k \gg 0$.  Let $F^{(m)}$ denote the linear
span of all semi-infinite monomials of charge $m$.  Then we have the
charge decomposition
\[
F = \bigoplus_{m \in \Z} F^{(m)}.
\]
We call $F^{(0)}$ \emph{fermionic Fock space}.

To any partition $\lambda = (\lambda_1 \ge \lambda_2 \ge \dots \ge
0)$ we associate a semi-infinite monomial $\varphi_\lambda =
\underline{i_0} \wedge \underline{i_1} \wedge \dots$ of charge $m$,
by letting $i_k = (m-k) + \lambda_k$.  This gives a bijection
between the set of all semi-infinite monomials of a fixed charge $m$
and the set $\P$ of all partitions (finite non-increasing sequences
of non-negative integers).  We define the \emph{energy} of
$\varphi_\lambda$ to be $|\lambda| := \sum_i \lambda_i$, the size of
the partition $\lambda$.  Let $F_j^{(m)}$ denote the linear span of
all semi-infinite monomials of charge $m$ and energy $j$. We then
have the energy decomposition
\[
F^{(m)} = \sum_{j \in \Z} F_j^{(m)}.
\]

For $j \in \Z$, define the \emph{wedging} and \emph{contracting}
operators $\psi_j$ and $\psi_j^*$ on $F$ by:
\begin{align*}
\psi_j(\underline{i_0} \wedge \underline{i_1} \wedge \dots) &=
\begin{cases}
0 & \text{if $j = i_s$ for some $s$,} \\
(-1)^{s+1} \underline{i_0} \wedge \dots \wedge \underline{i_s}
\wedge \underline{j} \wedge \underline{i_{s+1}} \wedge \dots &
\text{if $i_s > j > i_{s+1}$}.
\end{cases} \\
\psi_j^*(\underline{i_0} \wedge \underline{i_1} \wedge \dots) &=
\begin{cases}
0 & \text{if $j \ne i_s$ for all $s$,} \\
(-1)^s \underline{i_0} \wedge \underline{i_1} \wedge \dots \wedge
\underline{i_{s-1}} \wedge \underline{i_{s+1}} \wedge \dots &
\text{if $j=i_s$}.
\end{cases}
\end{align*}
These can be thought of as creation and annihilation operators. The
operator $\psi_j$ creates a particle in state $j$ while the operator
$\psi_j^*$ annihilates a particle in state $j$. Note that
\[
\psi_j (F^{(m)}) \subseteq F^{(m+1)},\quad \psi_j^*(F^{(m)})
\subseteq F^{(m-1)}.
\]
The operators $\psi_j$ and $\psi_j^*$ are called \emph{free
fermions}.  One can check directly that $\psi_j$ and $\psi_j^*$ are
adjoint with respect to the Hermitian form $H(\cdot,\cdot)$ and that
the following relations hold:
\begin{equation}
\label{eq:cliff-relations} \psi_i \psi_j^* + \psi_j^* \psi_i =
\delta_{ij},\ \psi_i \psi_j + \psi_j \psi_i = 0,\ \psi_i^* \psi_j^*
+ \psi_j^* \psi_i^* = 0.
\end{equation}
Thus, the operators $\psi_j$ and $\psi_j^*$ generate a Clifford
algebra $\Cl$.  It is easily seen that $F$ is an irreducible
$\Cl$-module and that
\[
\psi_j \vac{m} = 0 \text{ for $j \le m$},\quad \psi_j^* \vac{m} = 0 \text{ for $j > m$}.
\]

Let $gl_\infty$ denote the Lie algebra of all complex (infinite)
matrices $(a_{ij})_{i,j \in \Z}$ such that the number of nonzero
$a_{ij}$ is finite, with the usual commutator bracket, and let
\[
sl_\infty = \{a \in gl_\infty\ |\ \tr a = 0\}.
\]
Let $E_{ij} \in gl_\infty$ denote the matrix with $(i,j)$-entry
equal to one and all other entries equal to zero.

There is an embedding $gl_\infty \to \Cl$ defined by
\[
r(E_{ij}) = \psi_i \psi_j^*
\]
and this defines a representation $r$ of $gl_\infty$ on $F$.  It is
easy to see that each $F^{(m)}$ is stable under the action of
$gl_\infty$ and thus $r$ restricts to a representation $r_m$ of
$gl_\infty$ on $F^{(m)}$ for each $m \in \Z$.  Note that $r(^t \bar
a)$ and $r(a)$, for $a \in gl_\infty$, are adjoint operators with
respect to the Hermitian form $H$.

One can check that $gl_\infty$ acts by derivations on $F$.  That is,
for $a = (a_{ij}) \in gl_\infty$,
\begin{equation}
\label{eq:gl-action} r(a)(\underline{i_0} \wedge \underline{i_1}
\wedge \dots) = (a \cdot \underline{i_0}) \wedge \underline{i_1}
\wedge \dots + \underline{i_0} \wedge (a \cdot \underline{i_1})
\wedge \dots,
\end{equation}
where we view $\underline j$ as the vector with $j$th component one
and all other components zero.  Thus $a \cdot \underline{j} = \sum_i
a_{ij} \underline{i}$.  One can then use the usual rules of the
exterior algebra to express the right-hand side of
\eqref{eq:gl-action} in terms of semi-infinite monomials.  Thus $F$
is an infinite generalization of the usual exterior algebra.  For
this reason, $r$ is often called the \emph{infinite wedge
representation}.

The representations $r_m$ are irreducible.  Also,
\begin{align*}
r_m(E_{ij}) \vac{m} &= 0, \text{ if $i<j$, or $i=j > m$},\\
r_m(E_{jj}) \vac{m} &= \vac{m}, \text{ if $j \le m$}.
\end{align*}
Thus, as an $sl_\infty$-module, $F^{(m)}$ is isomorphic to the
irreducible integral representation $L(\omega_m)$ of highest weight
$\omega_m$, where $\omega_m$ is the fundamental weight defined by
$\left< \omega_m, \alpha_j^\vee \right> = \delta_{mj}$.  Here the
$\alpha_j^\vee$ are the simple coroots and $\left< \cdot, \cdot
\right>$ is the usual pairing between weights and coweights.  The
representation $r_0$ of $gl_\infty$ (or $sl_\infty$) on $F^{(0)}$ is
called the \emph{basic representation}.

We can define an $\mathfrak{s}$-module structure on $F^{(m)}$ by
introducing the \emph{free bosons}
\begin{align*}
\alpha_n &= \sum_{j \in \Z} \psi_j \psi_{j+n}^*,\ n \in \Z
\backslash \{0\}, \\
\alpha_0 &= \sum_{j > 0} \psi_j \psi_j^* - \sum_{j \le 0} \psi_j^*
\psi_j.
\end{align*}
Note that while the sums involved in the above definitions are
infinite, all but a finite number of them act as zero on any
semi-infinite monomial and so the operations are well-defined.

\begin{prop}[{\cite[Proposition~14.9]{K}}] \label{prop:fermionic-smod}
The map $s_n \mapsto \alpha_n$, $n \in \Z$, $K \mapsto \Id$ defines
an $\mathfrak{s}$-module structure on $F^{(m)}$ and $F^{(m)}$ is
irreducible for all $m \in \Z$.
\end{prop}

Note that
\[
\alpha_0 |_{F^{(m)}} = m\id
\]
and that $\alpha_m$ and $\alpha_{-m}$ are adjoint operators.

It follows from Proposition~\ref{prop:fermionic-smod} and a
uniqueness property for representations of $\mathfrak{s}$ (see
\cite[Corollary~9.13]{K}) that there is a unique isomorphism of
$\mathfrak{s}$-modules
\[
\sigma : F \stackrel{\cong}{\rightarrow} B,
\]
such that $\sigma(\vac{m}) = q^m$.  Note that
\[
\sigma(F^{(m)}) = B^{(m)} := q^m \C [p_1,p_2,\dots].
\]
We will denote this restriction of $\sigma$ to $F^{(m)}$ by
$\sigma_m$.  Also, we call $B^{(0)}$ \emph{bosonic Fock space}.

For a partition $\lambda \in \P$, we recall the definition of the
\emph{Schur polynomial} $S_\lambda$.  First, one defines the
elementary Schur polynomials $S_n$ by
\begin{gather*}
S_n = 0 \text{ for } n<0,\quad S_0 = 1, \\
S_n = \sum_{\mu \vdash n} \frac{1}{z_\lambda} p_1^{m_1(\mu)}
p_2^{m_2(\mu)} \dots \text{ for } n > 0.
\end{gather*}
Where $\mu \vdash n$ means that $\mu$ is a partition of size $n$.
Recall that $m_i(\mu)$ is the number of parts of $\mu$ equal to $i$.
Then for $\lambda = (\lambda_1 \le \lambda_2 \le \dots) \in \P$
define
\[
S_\lambda = \det(S_{\lambda_i + j - i})_{1 \le i,j \le |\lambda|}.
\]
Note that to translate to the definition of Schur polynomials in the
ring of symmetric functions in variables $x_1, x_2, \dots$ (see
\cite{Mac95}), one should replace $p_i$ with the $i$th power sum.

If $\varphi_\lambda \in F^{(m)}$ is a semi-infinite monomial, then
(see \cite[Theorem~14.10]{K})
\[
\sigma(\varphi_\lambda) = q^m S_{\lambda}.
\]

The process we have just outlined, constructing bosons in terms of
fermions acting on full fermionic Fock space, is called
\emph{bosonization}.  There exists an opposite procedure,
\emph{fermionization}, which consists of constructing fermions in
terms of bosons acting on full bosonic Fock space.  This task is
somewhat more complicated, involving vertex operator algebras. Since
we do not need fermionization in the current paper, we will not
describe the procedure here but instead refer the reader to
\cite[{\S\S\ 14.9-14.10}]{K}.


\section{Equivariant cohomology and localization}
\label{sec:localization}

In this section we recall the definition of equivariant cohomology
and the localization theorem.  This will be our main tool in
relating the geometric realizations of the bosonic and fermionic
Fock spaces thus yielding a geometric boson-fermion correspondence.
We concentrate on the case where the group is the torus and follow
the presentation in \cite{CK}.

Let $T=\C^*$ be the one-dimensional torus.  For an algebraic variety
$X$ equipped with a $T$-action, let $H^*_T(X)$ denote the
equivariant cohomology ring of $X$ with complex coefficients.  We
recall the definition. Let $BT = \CP^\infty$ and $ET$ be the
tautological bundle on $\CP^\infty$. Set $X_T = X \times_T ET$. This
is a bundle over $BT$ with fiber $X$. Then, by definition,
\[
H_T^*(X) = H^*(X_T),
\]
where $H^*(X_T)$ is the ordinary cohomology of $X_T$.

Recall that we have flat equivariant pullbacks and proper
equivariant pushforwards in equivariant cohomology.  If $\text{pt}$
is the space consisting of a single point with the trivial
$T$-action, then $H_T^*(\text{pt}) = \C[t]$ where $t$ is an element
of degree 2. Thus, by the pullback via $M \to \text{pt}$, $H^*_T(X)$
has the structure of a $\C[t]$-module.  For a proper $T$-equivariant
morphism $f : Y \to X$ of algebraic varieties, we have a Gysin map
$f_! : H_T^*(Y) \to H_T^*(X)$.  If $Y$ is a $T$-equivariant
codimension-$k$ closed subvariety of $X$ and $i:Y \hookrightarrow X$
is the inclusion map, we define
\[
[Y] = i_!(1_Y) \in H^{2k}_T(X),
\]
where $1_Y \in H_T^0(Y)$ is the unit in $H_T^*(Y)$

An \emph{equivariant vector bundle} is a vector bundle $E$ over $X$
such that the action of $T$ on $X$ lifts to an action of $E$ which
is linear on the fibers.  Then $E_T$ is a vector bundle over $X_T$
and the \emph{equivariant Chern classes} $c^T_k(E) \in H_T^*(X)$ are
defined to be the ordinary Chern classes $c_k(E_T)$.  If $E$ has
rank $r$, then the top Chern class $c_r^T(E)$ is called the
\emph{equivariant Euler class} of $E$ and is denoted $e_T(E) \in
H_T^*(X)$.

Very important in our discussion will be the \emph{localization
theorem} which we now describe.  Suppose that $X$ is smooth and has
a $T$-action.  Then the fixed point locus $X^T$ is a union of smooth
connected components $Z_j$.  Let $i_j : Z_j \hookrightarrow X$ be
the inclusion and let $N_j$ denote the normal bundle of $Z_j$ in
$X$.  Then $N_j$ is an equivariant vector bundle and thus has an
equivariant Euler class
\[
e_T(N_j) \in H_T^*(Z_j).
\]
The equivariant inclusion $i_j : Z_j \to X$ induces the pullback map
$i_j^* : H_T^*(X) \to H_T^*(Z_j)$.  The Gysin map
\[
i_{j!} : H_T^*(Z_j) \to H_T^*(X)
\]
has the property that for any $\alpha \in H_T^*(Z_j)$,
\[
i_j^* \circ i_{j!} (\alpha) = \alpha \cup e_T(N_j).
\]
Let $\C(t)$ be the field of fractions of $\C[t]$ and form the
\emph{localization} $H_T^*(Z_j) \otimes \C(t)$.  Then $e_T(N_j)$ is
invertible in $H_T^*(Z_j) \otimes \C(t)$.  The following proposition
is referred to as the \emph{localization theorem}.

\begin{prop}[{\cite[Proposition~9.1.2]{CK},\cite{AB}}]
There is an isomorphism
\[
H_T^*(X) \otimes \C(t) \stackrel{\simeq}{\longrightarrow}
\bigoplus_j H_T^*(Z_j) \otimes \C(t)
\]
given by
\[
\alpha \mapsto (i_j^*(\alpha)/e_T(N_j))_j.
\]
The inverse map is given by
\[
(\alpha_j)_j \mapsto \sum_j i_{j!}(\alpha_j).
\]
In particular, for any $\alpha \in H_T^*(X) \otimes \C(t)$, we have
\[
\alpha = \sum_j i_{j!} \left( \frac{i_j^*(\alpha)}{e_T(N_j)}
\right).
\]
\end{prop}


\section{Geometric realization of fermionic Fock space}
\label{sec:geometric-fock-space}

In this section we describe a geometric realization of fermionic
Fock space using the quiver varieties of Nakajima and the results of
\cite{FS03}.  We only introduce here the special case of quiver
varieties corresponding to the basic representation $L(\omega_0)$ of
the Lie algebra $sl_\infty$ of type $A_\infty$.  In this case, the
quiver varieties are simply points.  However, the reader should keep
in mind the fact that the construction generalizes to irreducible
integrable representations of symmetric Kac-Moody algebras (see
\cite{N94,N98}).  Note that we use a different stability condition
that the one used in \cite{N94,N98} and so our definitions differ
slightly from the ones that appear there. One can translate between
the two stability conditions by taking transposes of the maps
appearing in the definitions of the quiver varieties. See \cite{N96}
for a discussion of various choices of stability condition.  Another
difference in our presentation below is that we use equivariant
cohomology rather than ordinary cohomology (or Borel-Moore
homology).  Since the varieties involved are points, this change is
minor.  However, we will need this formulation to connect the quiver
variety picture to the geometric construction of the bosonic Fock
space described later.

Let $V = \bigoplus_{k \in \Z} V_k$ be a finite dimensional complex
$\Z$-graded vector space of graded dimension $\v = (\dim V_k)_{k \in
\Z}$. Then we define $\mathbf{M}(\v)$ to be the set of all triples
$(B_1,B_{-1},\tv)$ where $\tv \in V_0$ and $B_1$ and $B_{-1}$ are
endomorphisms of the graded vector space $V$ of degrees 1 and $-1$
respectively satisfying
\[
[B_1,B_{-1}] = B_1B_{-1} - B_{-1} B_1 = 0.
\]

Now, let
\[
G_\v = \prod_{k \in I} GL(V_k)
\]
be the group of grading-preserving automorphisms of $V$.  Then we
have a natural action of $G_\v$ on $\mathbf{M}(\v)$ given by
\[
g \cdot (B_1,B_{-1},\tv) = (g B_1 g^{-1}, g B_{-1} g^{-1}, g(\tv)).
\]

We say that a graded subspace $S$ of $V$ is
\emph{$(B_1,B_{-1})$-invariant} if $B_1(S) \subset S$ and $B_{-1}(S)
\subset S$.  We say that a point $(B_1, B_{-1},\tv)$ of
$\mathbf{M}(\v)$ is \emph{stable} if any $(B_1,B_{-1})$-invariant
graded subspace $S$ of $V$ containing $\tv$ is equal to all of $V$.

We let $\mathbf{M}(\v)^s$ denote the set of all stable points.  It
is known (see \cite[{Lemma~3.10}]{N98}) that the stabilizer in
$G_\v$ of any point in $\mathbf{M}(\v)^s$ is trivial and we define
the \emph{quiver variety}
\[
\M(\v) = \mathbf{M}(\v)^s/G_\v.
\]
This is a geometric quotient (it can also be viewed as a symplectic
quotient, although we do not discuss the symplectic structure here).
For $(B_1,B_{-1},\tv) \in \mathbf{M}(\v)^s$, we denote the
corresponding orbit in $\M(\v)$ by $[B_1,B_{-1},\tv]$.

\begin{figure}
\centering \epsfig{file=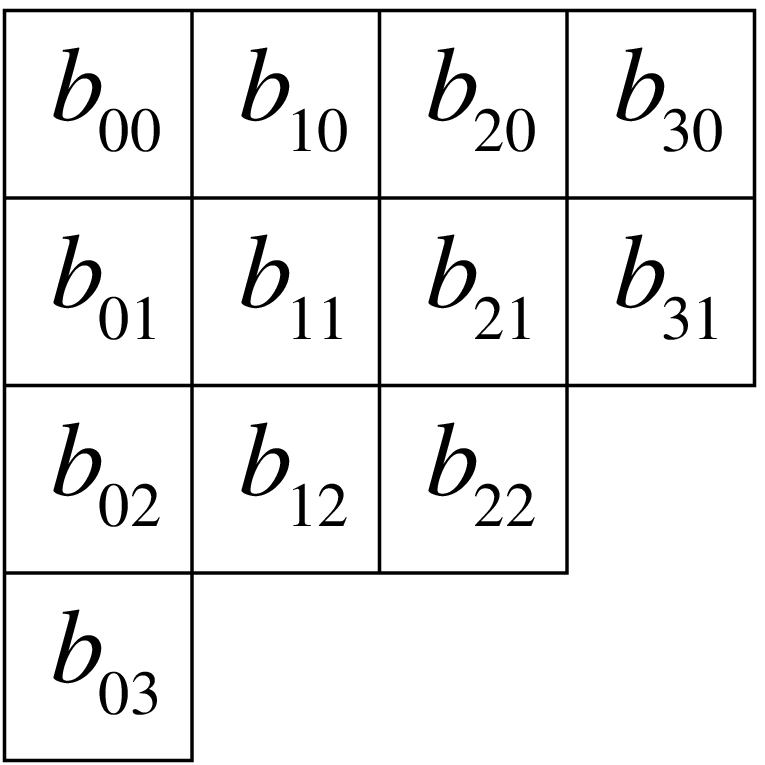,width=0.2\textwidth}
\caption{The labels of the boxes of the Young diagram
$D_{(4,4,3,1)}$. \label{fig:youngdiaglabels}}
\end{figure}

For a partition $\lambda \in \P$, let $D_\lambda$ denote the Young
diagram corresponding to $\lambda$.  We view $D_\lambda$ as a
collection of boxes, with $b_{jk} \in D_\lambda$ denoting the box in
the $j$th column and $k$th row where we start numbering from zero.
For example, the labels of the boxes of the Young diagram
$D_{(4,4,3,1)}$ are as in Figure~\ref{fig:youngdiaglabels}.

\begin{figure}
\begin{center}
\epsfig{file=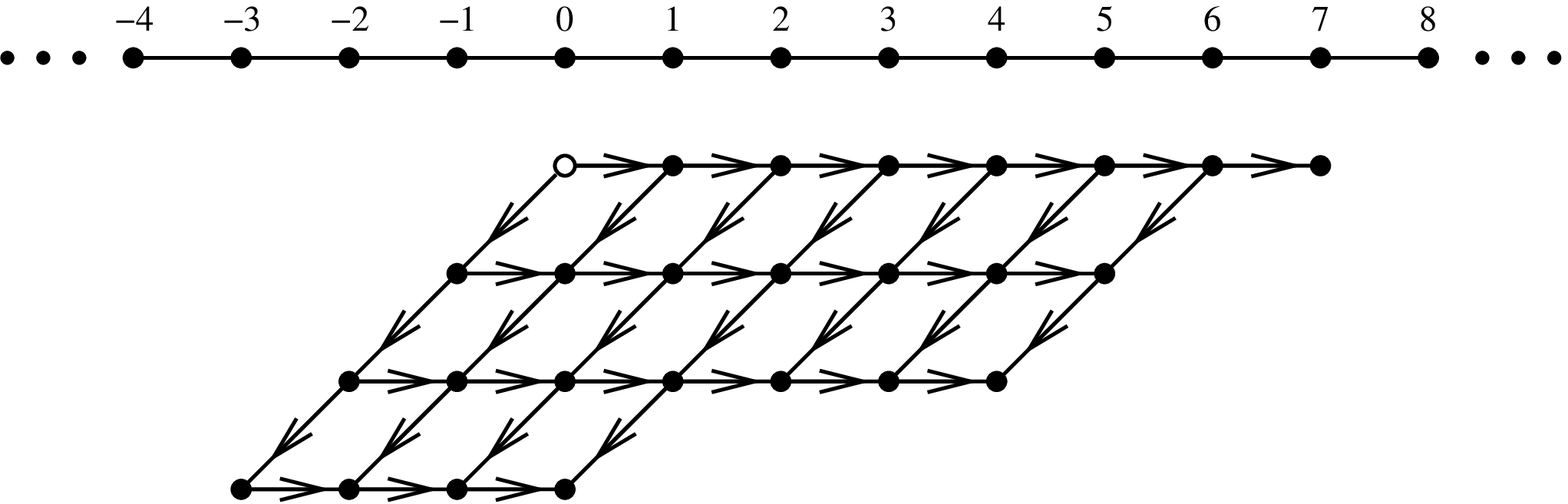,width=4in} \caption{A pictoral
representation of $(B_1^\lambda,B_{-1}^\lambda,\tv^\lambda)$ for
$\lambda=(8,7,7,4)$. The top line is the Dynkin graph of type
$A_\infty$. The vertices below represent the boxes $b_{jk} \in
D_\lambda$ while the arrows represent the actions of $B_1^\lambda$
and $B_{-1}^\lambda$. The vector space $V^\lambda_k$ is spanned by
the vertices directly beneath the vertex $k$ in the Dynkin diagram
and the vector $\tv^\lambda = b_{00}$ is indicated by a hollow
vertex.} \label{fig:quiver-youngdiag}
\end{center}
\end{figure}

We define the \emph{residue} of a box $b_{jk}$ to be $j-k$.  Define
$V^\lambda_k$ to be the $\C$-span of the boxes in $D_\lambda$ of
residue $k$.  Then $\dim V^\lambda = \v^\lambda$ where $v^\lambda_k$
is the number of boxes in $D_\lambda$ of residue $k$.  Define an
element $(B_1^\lambda,B_{-1}^\lambda,\tv^\lambda)$ of $M(\v)$ by
\begin{align*}
B^\lambda_1(b_{j,k}) &= b_{j+1,k},\ j,k \in \Z_{\ge 0} \\
B^\lambda_{-1}(b_{j,k}) &= b_{j,k+1},\ j,k \in \Z_{\ge 0} \\
\tv^\lambda &= b_{00},
\end{align*}
where $b_{jk}=0$ if $b_{jk} \not \in D_\lambda$.  We picture
$(B_1^\lambda,B_{-1}^\lambda,\tv^\lambda)$ as in
Figure~\ref{fig:quiver-youngdiag}.

In was shown in \cite{FS03} that $\M(\v)$ is empty unless $\v =
\v^\lambda$ for some $\lambda \in \P$.  Also, $\M(\v^\lambda) =
[B_1^\lambda,B_{-1}^\lambda,\tv^\lambda]$ is a single point.  This
follows from the fact that, in general, the quiver variety $\M(\v)$
associated to a Kac-Moody algebra with symmetric Cartan matrix $C$
is connected (see \cite[Theorem~6.2]{N98}) and its dimension is
given by (see \cite[Corollary~3.12]{N98})
\[
\dim_\C \M(\v) = 2v_0 - \v \cdot C \v.
\]
For $C$ the Cartan matrix of type $A_\infty$ (i.e. $C_{ij} =
2\delta_{ij} - \delta_{i-1,j} - \delta_{i+1,j}$), one can see from
this formula that $\dim M(\v^\lambda) = 0$ for $\lambda \in \P$.
This is, of course, not true in general.  For other types, quiver
varieties can have higher dimension.

Let the torus $T=\C^*$ act trivially on the quiver varieties.  Since
the equivariant cohomology of a point $H^*_T(pt)$ is $\C[t]$, we
have the following.
\[
H^*_T(\M(\v)) \cong \begin{cases} \C[t]
& \text{if $\v = \v^\lambda$ for $\lambda \in \P$} \\
0 & \text{otherwise} \end{cases}.
\]

For $k \in I$, define the \emph{Hecke correspondence} $\B_k(\v)$ to
be the variety of all $(B_1,B_{-1},\tv,S)$ (modulo the
$G_\v$-action) such that $(B_1,B_{-1},\tv) \in \mathbf{M}(\v)^s$ and
$S$ is a $(B_1,B_{-1})$-invariant subspace such that $\dim S =
\mathbf{e}^k$ where $\mathbf{e}^k$ has $k$-component equal to one
and all other components equal to zero. We consider the $G_\v$-orbit
through $(B_1,B_{-1},\tv,S)$ as a point in $\M(\v-\e^k) \times
\M(\v)$ by factoring by the subspace $S$ in the first factor.  Note
that from the explicit description of the $M(\v)$ given above, we
see that $\B_k(\v)$ is empty unless $\v = \v^\lambda$ for some
$\lambda \in \P$ such that there exists a $\mu \in \P$ with $D_\mu$
differing from $D_\lambda$ by the removal of a single box of residue
$k$. Then
\[
\B_k(\v^\lambda) = ([B_1^\mu,B_{-1}^\mu,\tv^\mu],
[B_1^\lambda,B_{-1}^\lambda, \tv^\lambda]).
\]

Let $\omega : \M(\v^1) \times \M(\v^2) \to \M(\v^2) \times \M(\v^1)$
be the map that interchanges the two factors.  We then define two
operators $E_k$ and $F_k$ that act on $\bigoplus_{\v} H^*_T(\M(\v))$
as follows.
\begin{gather}
E_k c = t^{-1} \cup {p_1}_! (p_2^* c \cap [\B_k(\v))],\quad c \in H^*_T (\M(\v)), \\
F_k c = t \cup {p_2}_![(p_1^*c \cap [\omega(\B_k(\v+\e^k,\w))]),
\quad c \in H^*_T (\M(\v)).
\end{gather}
Here we have used the projection maps
\[
\M(\v_1) \stackrel{p_1}{\longleftarrow} \M(\v_1) \times \M(\v_2)
\stackrel{p_2}{\longrightarrow} \M(\v_2).
\]
Note that the operators $E_k$ and $F_k$ preserve the subspace
\[
\H^F := \bigoplus_{\lambda \in \P} H_T^{2|\lambda|}(\M(\v^\lambda)).
\]
We also define
\[
\H^F_n = \bigoplus_{\lambda \vdash n} H_T^{2n}(\M(\v^\lambda)) =
H^{2n}_T \left( \bigsqcup_{\lambda \vdash n} \M(\v^\lambda) \right).
\]
Note that since the $\M(\v^\lambda)$ are all points, we have
\[
H^{2k}_T(\M(\v^\lambda)) = t^k \cup H_T^0(\M(\v^\lambda)),\ k \in
Z_{\ge 0},\ \lambda \in \P.
\]
In particular, we see from this that $E_k$ as defined above is
indeed an operator on $\H^F$ (a priori, it is an operator on $\H^F
\otimes \C(t)$).

\begin{theo}
\label{thm:nak-hom} The operators $E_k$ and $F_k$ satisfy the
relations of the Chevalley generators of $sl_\infty$ and thus define
an action of $sl_\infty$ on $\H^F$.  Under this action, $\H^F$ is
the basic representation. The class $[\M(\mathbf{0})]$ is a highest
weight vector and $H^{2|\lambda|}(\M(\v^\lambda))$ is the weight
space of weight $\omega_0 - \sum v^\lambda_k \alpha_k$, where
$\alpha_k$ are the simple roots of $sl_\infty$.
\end{theo}

\begin{proof}
This follows immediately from the results of \cite{N98}.  Our
modifications, as noted at the beginning of this section, are minor.
\end{proof}


\section{The Hilbert scheme, torus action and fixed points}
\label{sec:hilbert-scheme-fixed-pts}

In this section we introduce the Hilbert scheme of $n$ points in
$\C^2$ and recall some results regarding the fixed points of the
natural torus action.  We shall see that these fixed points are
naturally identified with the quiver varieties of
Section~\ref{sec:geometric-fock-space}.  This precise relationship
between the spaces involved in the geometric constructions of the
bosonic and fermionic Fock spaces will allow us to use the
localization theorem to yield a geometric boson-fermion
correspondence.  For other results relating quiver varieties and
Hilbert schemes, we refer the reader to \cite{Hai,Kuz01,VV99,Wan02}.

Consider the space of $n$ (unordered) points in $\C^2$.  This is a
singular space, the singularities occurring when points collide. The
Hilbert scheme of $n$ points in $\C^2$ is a resolution of
singularities of this space.  It resolves the singularities by
``retaining information about how points collided''.  It is thus a
very natural space and has been used to construct representations of
infinite dimensional Heisenberg algebras and Virasoro algebras. We
describe here the former.

Precisely, the Hilbert scheme $\hs$ parameterizes 0-dimensional
closed subschemes of $\C^2$ of length $n$. So if $x$ and $y$ are the
standard coordinate functions on $\C^2$,
\[
\hs = \{I\ |\ I \text{ is an ideal of } \C[x,y],\ \dim \C[x,y]/I =
n\}.
\]
Consider the action of the one-dimensional torus $T=\C^*$ on $\C^2$
by
\[
z \cdot (x,y) = (z x, z^{-1} y),\quad z \in T.
\]
The only fixed point is the origin, which we will denote by $u$.

The $T$-action on $\C^2$ induces a $T$-action on $\hs$.  The support
of a $T$-fixed point in $\hs$ is $u$ since it must be a fixed point
of $\C^2$.  In order to explicitly describe the $T$-fixed points and
relate them to the quiver varieties described above, we give the
following alternate description of the Hilbert schemes (see
\cite[Theorem~1.9]{Nak99}).
\begin{equation}
\label{eq:hs-description} \hs \cong \left\{(B^1,B^2,\tv)\ \left| \
[B^1,B^2] = 0,\ (B^1,B^2,\tv) \text{ is stable} \right.
\right\}/GL(V),
\end{equation}
where $B^j \in \End(V)$, $\tv \in V = \C^n$, and we say that
$(B^1,B^2,\tv)$ is stable if there exists no proper subspace $S
\varsubsetneq V$ such that $B^j(S) \subset S$ for $j = 1,2$, and
$\tv \in S$.  The action of $GL(V)$ is given by
\[
g \cdot (B^1,B^2,\tv) = (g B^1 g^{-1}, g B^2 g^{-1}, g(\tv)),\quad g
\in GL(V).
\]
In this description, the action of $T$ on $\hs$ is given by
\[
z \cdot [B^1,B^2,\tv] = [z B^1, z^{-1} B^2, \tv],\quad z \in T.
\]
Here $[B^1,B^2,\tv]$ denotes the $GL(V)$-orbit through
$(B^1,B^2,\tv)$.

From another description of the Hilbert scheme (see
\cite[Theorem~3.24]{Nak99}) one can see that $[B^1,B^2,\tv] \in \hs$
is a fixed point if and only if there exists a homomorphism $\lambda
: T \to U(\C^n)$ such that
\begin{align*}
z B^1 &= \lambda(z)^{-1} B^1 \lambda(z), \\
z^{-1} B^2 &= \lambda(z)^{-1} B^2 \lambda(z), \\
\tv &= \lambda(z)^{-1} (\tv).
\end{align*}
Thus, if $[B^1,B^2,\tv]$ is a fixed point, we have a weight
decomposition of $V$ with respect to $\lambda(z)$ given by
\[
V = \bigoplus_k V_k,
\]
where
\[
V_k = \{v \in V\ |\ \lambda(z) \cdot v = z^{-k} v\}.
\]
It follows that the only non-zero components of $B^1$ and $B^2$ are
\begin{align}
\label{eq:fpc1} B^1 &: V_k \to V_{k+1}, \\
\label{eq:fpc2} B^2 &: V_k \to V_{k-1}
\end{align}
and
\begin{equation}
\label{eq:fpc3} \tv \in V_0.
\end{equation}
So $B^1$ and $B^2$ are endomorphisms of the graded vector space $V$
of degrees 1 and $-1$ respectively.  Conversely, one can see that
any triple $(B^1,B^2,\tv)$ satisfying
conditions~\eqref{eq:fpc1}-\eqref{eq:fpc3} for some $\lambda$ is a
$T$-fixed point. This precisely matches the description of the
points of the $A_\infty$ quiver varieties given in
Section~\ref{sec:geometric-fock-space}.  We simply set $B_1 = B^1$
and $B_{-1} = B^2$ and then the $T$-fixed points of the Hilbert
scheme are exactly the quiver varieties of type $A_\infty$. In
particular, the fixed point set consists of isolated points.  Recall
that the total dimension of the vector space $V$ appearing in the
definition of the quiver variety $\mathfrak{M}(\v^\lambda)$ is
$|\lambda|$, the size of the partition $\lambda$.  Thus there is a
natural identification
\[
\hs^T = \bigsqcup_{\lambda \vdash n} \M(\v^\lambda)
\]
and hence a natural identification
\[
\bigoplus_n H^{2n}_T(\hs^T) = \H^F.
\]


\section{Geometric realization of bosonic Fock space}
\label{sec:geom-bosonic}

In this section we describe the geometric realization of bosonic
Fock space in the equivariant cohomology of the Hilbert schemes of
points in $\C^2$.  We refer the reader to
\cite{Gro96,LQW,Nak96b,Nak99,Vas01} for more details.

Let $\Sigma$ be the $x$-axis in $\C^2$ and recall that $u$ is the
origin.  As a $T$-module, we have $T_u\Sigma = \theta^{-1}$ where
the $\theta$ is the one-dimensional standard $T$-module. By the
localization theorem,
\begin{equation} \label{eq:sigma-class}
[\Sigma] = -t^{-1}[u].
\end{equation}

The tangent space of $\hs$ at the $T$-fixed point $\M(\v^\lambda)$
is $T$-equivariantly isomorphic to (see
\cite[Proposition~5.8]{Nak99} and \cite{ES,Nak96b})
\[
T_{\M(\v^\lambda)} \hs = \bigoplus_{b \in D_\lambda} \left(
\theta^{\text{hook}(b)} \oplus \theta^{-\text{hook}(b)} \right)
\]
where $\text{hook}(b)$ is the hook length of the box $b$. Recall
that the hook length of a box is the number of boxes directly to its
right plus the number of boxes directly below it plus one (for the
box itself). For example, suppose $b$ is the indicated box in the
following Young diagram.
\[
\young(\ \ \ \ \ \ ,\ \ b\bullet\bullet\bullet,\ \ \bullet\ \ ,\ \
\bullet,\ \ )
\]
Then $\text{hook}(b)=6$ (the number of boxes containing bullets plus
one for the box containing $b$). Thus
\begin{equation}
\label{eq:euler-class-fixed-pt} e_T(T_{\M(\v^\lambda)} \hs) = (-1)^n
h(\lambda)^2 t^{2n},
\end{equation}
where
\[
h(\lambda) = \prod_{b \in D_\lambda} \text{hook}(b).
\]
Note that $[\M(\v^\lambda)] \in H^{4n}_T(\hs)$.  Define
\[
[\lambda] = \frac{(-1)^n}{h(\lambda)} t^{-n} [\M(\v^\lambda)].
\]

The odd Betti numbers of $\hs$ are equal to zero and $H^k(\hs) = 0$
for $k > 2n$.  We have
\[
H^{2k}_T(\hs) = t^{k-n} \cup H^{2n}_T(\hs),\ k \ge n.
\]
Define
\[
\H^B_n = H^{2n}_T(\hs),\quad \H^B = \bigoplus_{n=0}^\infty \H^B_n.
\]

We now introduce a bilinear form on the equivariant cohomology.  Let
\[
\iota : \hs^T \hookrightarrow \hs
\]
denote the inclusion of the torus fixed points.  Then we have the
Gysin map
\[
\iota_! : H^*_T(\hs^T) \to H^*_T(\hs).
\]
We denote the induced map $H^*_T(\hs^T)' \to H^*_T(\hs)'$ also by
$\iota_!$ where
\begin{align*}
H^*_T(\hs^T)' &= H^*_T(\hs^T) \otimes \C(t),\\
H^*_T(\hs)' &= H^*_T(\hs) \otimes \C(t)
\end{align*}
are the localizations.  This is an isomorphism by the localization
theorem.

Define the bilinear form $\left< \cdot, \cdot \right> : H^*_T(\hs)'
\times H^*_T(\hs)' \to \C(t)$ by
\[
\left< \alpha, \beta \right> = (-1)^n p_! \iota_!^{-1} (\alpha \cup
\beta)
\]
where $p : \hs^T \to \text{pt}$ is the projection to a point. This
induces a bilinear form on
\[
\H' = \bigoplus_{n=0}^\infty H^*_T(\hs)'.
\]
For $\lambda, \mu \vdash n$, we have that (see \cite[Lemme~2]{Vas01}
and \cite[Equation~2.21]{LQW})
\[
\left< [\lambda], [\mu] \right> = \delta_{\lambda,\mu}.
\]
Now, by the localization theorem, the classes $[\lambda]$, $\lambda
\vdash n$, form a linear basis of the $\H^B_n$.  Thus, the
restriction to $\H^B_n$ of the bilinear form $\left< \cdot, \cdot
\right>$ is non-degenerate.  This induces a non-degenerate bilinear
form $\left< \cdot, \cdot \right> : \H^B \times \H^B \to \C$.

We now discuss the action of the oscillator algebra $\mathfrak{s}$
on the cohomology of the Hilbert schemes.  Define
\[
\Sigma_{n,i} = \{(I_1,I_2) \in \hs[n+i] \times \hs\ |\ I_1 \subset
I_2,\, \text{Supp} (I_2/I_1) = \{z\},\, z \in \Sigma\}.
\]
We have the two natural projections
\[
\hs[n+i] \stackrel{\pi_1}{\longleftarrow} \hs[n+i] \times \hs
\stackrel{\pi_2}{\longrightarrow} \hs.
\]
The restriction of $\pi_1$ to $\Sigma_{n,i}$ is proper.  Thus we can
form the linear operator $\p_{-i} \in \End \H'$ (see \cite{Vas01})
by
\[
\p_{-i}(\alpha) = \pi_{1!}(\pi_2^* \alpha \cup [\Sigma_{n,i}]),\quad
\alpha \in H^*_T(\hs)'.
\]
We then define $\p_i \in \End \H'$ to be the adjoint operator to
$\p_{-i}$ with respect to the bilinear form $\left< \cdot, \cdot
\right>$ on $\H'$.  One can show that
\[
\p_i(\alpha) = (-1)^i \pi'_{2!} (\iota \times \text{Id})_!^{-1}
(\pi_1^* \alpha \cup [\Sigma_{n-i,i}]),\quad \alpha \in H^*_T(\hs)',
\]
where $\pi'_2$ is the natural projection $\hs^T \times \hs[n-i] \to
\hs[n-i]$.  We also set $\p_0 = 0$.

For $i > 0$, the restriction $\p_{-i}$ to $\H^B$ yields a linear
operator in $\End \H^B$.  The restriction of $\p_i$ to $\H$ is the
adjoint operator to $\p_{-i}$ with respect to the non-degenerate
bilinear form $\left< \cdot, \cdot \right> : \H^B \otimes_\C \H^B
\to \C$, which is the restriction of the bilinear form $\left<
\cdot, \cdot \right>$ on $\H'$. Thus the restriction of $\p_i$ to
$\H^B$ is an operator in $\End \H^B$, again denoted by $\p_i$.

\begin{prop}[{\cite[Lemme~1]{Vas01}}] \label{prop:heis-action}
The operators $\p_k$, $k \in \Z$, acting on $\H^B$ satisfy the
following Heisenberg commutation relations:
\[
[\p_k,\p_l] = k \delta_{k,-l} \text{Id}.
\]
In particular,
\begin{align*}
s_k &\mapsto \p_k,\ k \in \Z,\\
K &\mapsto \text{Id},
\end{align*}
defines an action of the oscillator algebra $\mathfrak{s}$ on
$\H^B$. Moreover, $\H^B$ is isomorphic to the bosonic Fock space
$B^{(0)}$. The unit $1 \in H^0_T(\hs[0])$ of $H^*_T(\hs[0])$
corresponds to the highest weight vector $1 \in B^{(0)}$.
\end{prop}


\section{A geometric boson-fermion correspondence}
\label{sec:geom-correspondence}

We now have all the tools necessary to describe the geometric
construction of the boson-fermion correspondence.  We have seen how
the quiver varieties used in the geometric realization of fermionic
Fock space can be naturally viewed as the set of torus fixed points
of the Hilbert scheme used in the geometric realization of the
bosonic Fock space.  We are then in a position to invoke the
localization theorem which gives an isomorphism between the
equivariant cohomology of the Hilbert scheme and its fixed point set
thus yielding our geometric version of the boson-fermion
correspondence.

Define a map $\eta : \H^F \to \H^B$ by
\[
\eta(\alpha) = \frac{(-1)^{|\lambda|}}{h(\lambda)}t^{-|\lambda|}
{i_\lambda}_!(\alpha) \in \H^B_{|\lambda|},\quad \alpha \in
H_T^{2|\lambda|}(\M(\v^\lambda)).
\]
The factor of $(-1)^{|\lambda|}t^{-|\lambda|}/h(\lambda)$, which
arises from the equivariant Euler class at the fixed point, ensures
that $\eta$ is an isometry if we endow $\H^F$ with the bilinear form
for which the $t^{|\lambda|} \cup 1_{\M(\v^\lambda)}$ form an
orthonormal basis.  It then follows from the localization theorem
and \eqref{eq:euler-class-fixed-pt} that $\eta$ is an isomorphism
with inverse given by
\[
\beta \mapsto \left(\frac{1}{h(\lambda)}t^{-n} i_\lambda^*(\beta))
\right)_{\lambda \vdash n} \in \H^F_n,\quad \beta \in \H^B_n.
\]

For a partition $\lambda = (\lambda_1 \ge \lambda_2 \ge \dots \ge
\lambda_l)$, define
\begin{align*}
z_\lambda &= \prod_{i \ge 1} i^{m_i(\lambda)} m_i(\lambda)!,\\
\p^{\lambda} &= \prod_{i \ge 1} \p_{-i}^{m_i(\lambda)},\\
\p_\lambda &= \p^\lambda \cdot 1 \in \H^B_n,\ \lambda \vdash n
\end{align*}
It follows from Proposition~\ref{prop:heis-action} and the fact that
$\p_{-i}$ and $\p_i$ are adjoint operators that
\[
\left< \p_\lambda, \p_\mu \right> = \delta_{\lambda,\mu} z_\lambda.
\]

\begin{theo}
\begin{enumerate}
\item \label{theo-bos} There exists an isomorphism $\phi : \H^B \to
B^{(0)}$ of $\mathfrak{s}$-modules preserving bilinear forms such
that
\[
\phi(\p_\lambda) = p_\lambda,\quad \phi([\lambda]) = S_\lambda.
\]

\item \label{theo-ferm}
There exists an isomorphism of $sl_\infty$-modules $\tau : F^{(0)}
\to \H^F$ such that
\[
\tau(\varphi_\lambda) = t^{|\lambda|} \cup 1_{\M(\v^\lambda)}.
\]

\item \label{theo-comp} We have $\phi \circ \eta \circ \tau = \sigma_0$.
\end{enumerate}
\end{theo}

\begin{proof}
Part~(\ref{theo-bos}) is proven in \cite[Proposition~2]{Vas01}.
Part~(\ref{theo-ferm}) follows from a comparison of the explicit
action $sl_\infty$ on the indicated bases in both spaces (see
\cite[Proposition~5.3]{FS03}).  Part~(\ref{theo-comp}) follows from
the fact that both $\phi \circ \eta \circ \tau$ and $\sigma_0$ are
linear isomorphisms sending a semi-infinite monomial
$\varphi_\lambda$ to $S_\lambda$.
\end{proof}

We can then define geometric bosons on the geometric fermionic Fock
space as follows.  For $\alpha \in \H^F$,
\[
s_k (\alpha) = \eta^{-1} \circ \p_k \circ \eta(\alpha),\quad K =
\Id.
\]
This defines a representation of $\mathfrak{s}$ on $\H^F$ and $\H^F
\cong B^{(0)}$ as $\mathfrak{s}$-modules.

We can also define an action of $sl_\infty$ on the geometric bosonic
Fock space by a similar procedure.  For $\beta \in \H^F$ and $k \in
\Z$,
\begin{align*}
E_k &= \eta \circ E_k \circ \eta^{-1}(\beta),\\
F_k &= \eta \circ F_k \circ \eta^{-1}(\beta).
\end{align*}

We thus see that in terms of the geometry of the Hilbert scheme, the
bosons correspondence to global operators while the fermions
correspond to local operators.

We should also note that the energy decomposition of the fermionic
Fock space has a nice geometric interpretation.  It corresponds to
grouping the quiver varieties according to the Hilbert scheme of
which they are a fixed point. More precisely,
\[
\tau (F^{(0)}_n) = \H^F_n = H^{2n}_T(\hs^T).
\]

Note that we have described the geometric boson-fermion
correspondence as a relationship between the basic representation of
$sl_\infty$ on fermionic Fock space $F^{(0)}$ and an irreducible
representation of the Heisenberg (or oscillator) algebra.  We could
also define the geometric action of the fermions themselves on the
full fermionic Fock space $F$. To do this we would have to introduce
the quiver varieties corresponding to the irreducible
representations $L(\omega_k)$ for $k \in \Z$. These are identical to
those for $L(\omega_0)$ except that one shifts the grading on the
vector space $V$. Then the fermions would be operators between the
cohomology of these different quiver varieties (see \cite{S03} for a
similar construction of Clifford algebras). However, the $sl_\infty$
action seems to be simpler geometrically and so we have emphasized
its role.

We should also comment on why the particular $T$-action we have
considered was chosen.  There is a natural action of the
two-dimensional torus $T^2 = (\C^*)^2$ on $\C^2$ given by
\[
(z_1,z_2) \cdot (x,y) = (z_1 x, z_2 y),\quad (z_1,z_2) \in T^2
\]
and this induces a $T^2$ action on the Hilbert scheme $X_n$.  The
$T$-action that we have considered arises from the embedding $T
\hookrightarrow T^2$ given by $z \mapsto (z,z^{-1})$.  A more
general embedding of $T$ into $T^2$ can certainly be considered.
This was examined in \cite{LQW,Nak96b} and one finds that the fixed
point contributions to the equivariant cohomology correspond to Jack
polynomials which generalize the Schur polynomials above.
Furthermore, one could consider the entire $T^2$ action and one
would obtain polynomials with an enlarged coefficient ring
containing a parameter which specialize to the Jack polynomials.
However, in these more general settings, one loses the natural
identification of the quiver varieties of type $A_\infty$ with the
set of fixed points. Under the embedding $z \mapsto (z,z^{-1})$, the
grading on the vector space $V$ in the definition of the quiver
variety appears naturally as described in
Section~\ref{sec:hilbert-scheme-fixed-pts}.


\bibliographystyle{abbrv}
\bibliography{biblist}

\end{document}